\documentclass[12pt,leqno]{article}\usepackage{amsmath}
\usepackage{amssymb}
\usepackage{amsmath}

\newcommand{\comment}[1]{}

\newcommand{\qed}{QED}

\newcommand{\Rng}{{\sf Rng}}

\newtheorem{lem}{Lemma}

\begin{document}
	
\title{A note on the submodel preservation property\\ in fragments of first-order logic}
	\author{Andr\'eka, H., van Benthem, J.\ and N\'emeti, I.}
	\date{}
	\maketitle

This  note contains some material promised in \cite{AvBNSub95, AvBNIGPL95, AvBNJPhL98}, along with some information on the current status of the problems mentioned in these papers. 
Section 1 contains an early example of failure of Los-Tarski for finite-variable fragments with binary relations from a 1992 manuscript, mentioned in \cite[p.10]{AvBNSub95} and in \cite[p.701]{AvBNIGPL95}, with no substantial change of content. For further background to this example, Section 2 surveys a number of results concerning the submodel preservation property for various fragments of first-order logic.

\section{A formula with three variables preserved under substructures but not equivalent to a universal formula}
In our first-order language we take one binary relation symbol, $R$.
Let  $(P,Z)$  denote the structure we will call the ``pentagon", namely

\begin{description}
\item{} $P=\{ 0,1,2,3,4\}$  and\
\item{} $Z=\{(i,j) : i<5\mbox{  and  }j=i+1(mod \,5)\mbox{ or }j=i-1(mod \,5)\}$.
\end{description}

Let  $(Q,T)$  denote the preceding structure but  5  replaced with  4,  i.e.,

\begin{description}
\item{} $Q=\{0,1,2,3\}$  and
\item{} $T=\{(i,j) : i<4\mbox{  and  }j=i+1(mod\, 4)\mbox{ or }j=i-1(mod\, 4)\}$.
\end{description}

We are going to define a first-order sentence $\varphi$  with 3 variables  defining (up to isomorphism) ``(P,Z) or one of its substructures", with a formulation given in Lemma \ref{1-lem}. Then  $\varphi$  is preserved by substructures,  $\varphi$  is true in  $(P,Z)$  and  $\varphi$  is not true in  $(Q,T)$. Next we will prove that every universal formula with three variables, constructed using atomic formulas and their negations, conjunction, disjunction and universal quantifiers,  which is true in  $(P,Z)$  is also true in $(Q,T)$, see Lemma \ref{2-lem}. Together, the  lemmas show that  $\varphi$  is not equivalent to any universal first-order formula with three variables.

\begin{lem}\label{1-lem} There is a sentence  $\varphi$  using only three variables such that in any structure  (W,S),
$\varphi$  is true in  (W,S)  if and only if
(W,S) is isomorphic to a substructure of  (P,Z).
\end{lem}
\noindent{\bf Proof.}  (1) \, Let  $Sxy$  denote  ``not $Rxy$ and not $x=y$".
Then define  $\sigma$  and  $\pi$  to be the following formulas:

\vspace{-1mm}

\begin{description}
\item \quad $\sigma = \forall xy\, [(Rxy\mbox{ iff }Ryx)\mbox{ and not }Rxx]$,
\item \quad $\pi  =  ``(R \circ R = S \cup Id)\mbox{ and }(R \circ S = R \cup S)\mbox{ and }(S \circ S = R \cup Id)"$,
\end{description}

\noindent where, e.g.,\  ``$(R \circ R = S \cup Id)$" abbreviates

\vspace{-1mm}

\[ \forall xy\, [\exists z(Rxz\mbox{ and }Rzy)\mbox{ iff }(Sxy\mbox{ or }x=y)]. \]

\vspace{1mm}

\noindent (Thus, ``$\circ$" abbreviates relation composition while ``$\cup$" abbreviates union.)

It is not difficult to check that $\sigma\land\pi$ is true only in one structure, namely, in the above-defined ``pentagon".

\vspace{2ex}
(2)\, Next, let  $q$  denote the following formula
\vspace{1ex}

``there exist x and y such that  [ \\

\vspace{-2ex}

 \quad \quad Sxy  and
 
  \quad \quad (there is only one z with Rxz) and
  
    \quad \quad(there is only one z with Ryz) and
    
  \quad \quad(there is a z such that Rxz and Syz)  and

  \quad \quad(there is a z such that Sxz and Ryz)  and
  
  \quad \quad(there is no z such that Sxz and Syz) and
  
  \quad \quad(there are x,y,z such that Rxy and Ryz)    ]".
  
\vspace{2ex}

It is not difficult to check that  $\sigma\land q$  is true only in one structure up to isomorphism, namely in the four-element substructure of the ``pentagon".

\vspace{3ex}

(3)\, Let  \mbox{tr}  denote the following formula

\medskip

``there exist x and y such that  [

\vspace{1mm}

\quad\quad Sxy  and

\quad\quad(there is only one z with Rxz)         and

\quad\quad(there is only one z with Ryz)         and

\quad\quad(there is a z such that Rxz and Ryz)   and

\quad\quad(there is no z such that Sxz and Syz) ]".
\medskip

\noindent Let  \mbox{ts}  denote the above formula but  R  and  S  interchanged everywhere. 

It is not difficult to check that the formula $\sigma\land(\mbox{tr}\lor\mbox{ts})$ is true only in the two three-element substructures of the ``pentagon".

\vspace{2ex}

(4)\, Let $\delta$  denote the following formula
\medskip

``there are at most two points".
\vspace{2ex}

(5) \,Finally, let  $\varphi$  denote the following  formula with three variables:

\vspace{-1mm}

\[ \varphi = \sigma \land (\pi \lor q \lor \mbox{tr} \lor \mbox{ts} \lor \delta) .\]

\vspace{1ex}

Now it can easily be checked that the  formula  $\varphi$  defined in this way has the properties required in Lemma \ref{1-lem}. \hfill \qed
\medskip

Next,  call a formula existential if it is built  from atomic formulas and their negations by using conjunction, disjunction, and existential quantifiers. Clearly, the negation of a universal formula in the above sense is equivalent to an existential formula containing the same number of variables.

\begin{lem}\label{2-lem} Every universal formula containing three variables and true in $(P,Z)$ is also true in $(Q,T)$.
\end{lem}

\noindent{\bf Proof.} Let  $V=\{ x,y,z\}$ denote the set of our variables, and let  $v:V \to  Q$ and  
$w:V\to P$  be two evaluations of the variables into $Q$ and $P$ respectively. We say that  $v$  and  $w$  are ``alike", if their ranges are isomorphic via the natural correspondence, i.e., if there is a bijection between  $\Rng v$  and  $\Rng w$  taking  $vx$, $vy$, and $vz$  to  $wx$, $wy$, and $w$z  respectively.

We will prove the following statement (*) by induction on  $\varphi$, where  $\varphi$ is any existential formula containing three variables:
\begin{description}
\item{(*)} $(Q,T) \models \varphi [v]$\,  implies\,  $(P,Z) \models \varphi [w]$,
for any alike evaluations  $v$  and  $w$.
\end{description}

If  $\varphi$  is an atomic formula or a negation of an atomic formula, then (*) is true because  $v$  and  $w$  are alike.  If  (*)  is true for the formulas  $\varphi$  and $\psi$,  then clearly  (*) is also true for  $\varphi\land\psi$  and  $\varphi\lor\psi$.

Now assume that (*) is true for  $\varphi$  and  let  $v$  and $w$  be alike evaluations. Suppose that  $(Q,T) \models\exists x \varphi [v]$.  This means that  
$(Q,T) \models \varphi [v(x\slash a)]$  for some element  $a$  of  $Q$,  where  $v(x\slash a)$  denotes the function we obtain from  $v$  by changing its value at  $x$  to $a$  and leaving the values of  $y$ and $z$  unchanged.  

Next, it is an easy observation that the two models  $(Q,T)$ and $(P,Z)$  have the following property:

\begin{description}
\item{(**)} For any alike  $v$  and  $w$, and for any element  $a$  of  $Q$  there is an element  $b$  of  $P$  such that  $v(x\slash a)$  and  $w(x\slash b)$  are again alike.
\end{description}

\noindent So let  $b$  be such that  $v(x\slash a)$  and $ w(x\slash b)$  are alike. Then, by our induction hypothesis on  $\varphi$,  $(P,Z)\models\varphi [w(x\slash b)]$  because  $(Q,T)\models\varphi [v(x\slash a)]$. Thus  $(P,Z) \models\exists x\varphi[w]$, and we are done with showing that, if (*) is true for  $\varphi$,  then (*) is also true for $\exists x\varphi$. The proofs for $\exists y$ and $\exists z$  are completely analogous.
This finishes the proof of (*).

\vspace{1ex}

Let   $\varphi$ now be any universal formula containing three variables and suppose that  $\varphi$  is not true in $(Q,T)$. We want to show that $\varphi$ is not true in  $(P,Z)$  either.  Let  $v$  be an evaluation of the variables such that  $(Q,T) \models\lnot\varphi[v]$. It is easy to check that, for any evaluation  $v$  into  $(Q,T)$,  there is an alike evaluation  $w$  into  $(P,Z)$. Take such a  $w$  which is alike to  $v$.  Then, since  $\lnot\varphi$ is equivalent to an existential formula, by (*) we have that  $(P,Z) \models\lnot\varphi [w]$,  i.e., that $\varphi$  is not true in $(P,Z)$. \hfill \qed
\bigskip

Finally, we note that the use of equality is not essential in the above example. We can obtain a formula without using equality and with the same properties by doing the following.
Let  $E$  be a new binary relation symbol in our language, and let  $\varepsilon$  denote the following formula:
\bigskip

$\varepsilon$ = ``($E$ is symmetric and transitive with domain the whole universe) 

\vspace{1mm}

and $(E \circ R = R \circ E = R$  and  $E \circ S = S \circ E = S$)",
\bigskip

\noindent where Sxy  denotes the formula  ``not Rxy and not Exy". Let  $\psi$  be the formula we obtain from the formula  $\varphi$  used in the preceding by changing  $x=y$  everywhere to  $Exy$,  and similarly for  $x=z$, etc.  Now the formula
$\varepsilon\land\psi$
is preserved by substructures, is not equivalent to any universal formula with three variables and it does not use the equality.

\section{Some related results and problems}
The property of a logic that  a formula is preserved under taking submodels if and only if it is equivalent to a universal formula is called the ({\L}os-Tarski) submodel preservation property (SPP for short). It is often also called the existential preservation property since in many logics it is equivalent to a formula being preserved under extensions of models if and only if it is equivalent to an existential formula. First-order logic FO with only relation symbols has these properties. The SPP may depend on how we formulate submodels and universal formulas since notions of universal formula that coincide for FO may  differ for other languages. In this note, `universal formulas' are built from atomic formulas and negations of atomic formulas by the use of conjunction, negation, and universal quantifiers (the `generalized universal formulas' of \cite[Definition 1.1]{AvBNSub95}). In FO, a formula is equivalent to a generalized universal formula if and only if only universal quantifiers occur in its prenex normal form. However, the two notions differ for the $n$-variable fragments of FO. 

We now list some results from the literature that may not be common knowledge,  to provide some background to the example in Section 1.

\vspace{1ex}

First-order logic with $n$ variables FO(n) does not have the SPP when $n\ge 2$. For $n\ge 3$ this is proved as \cite[Theorem 3.1]{AvBNSub95} (see also \cite[Theorem 3.2]{AvBNIGPL95}, \cite[Theorem 3.3.1]{AvBNJPhL98}). Failure of SPP for FO(2) is proved in \cite{GR99} as a solution to \cite[Problem 4.1]{AvBNSub95} (for this problem, see also \cite[item (1) on p.701]{AvBNIGPL95}, \cite[item (1) on p.241]{AvBNJPhL98}). For FO(1), it is easy to see that SPP holds with universal formulas understood as above but fails with universal formulas understood as having universal quantifiers only in a prenex normal form. 

For showing failure of SPP, the question of what kinds of relation symbols are needed is relevant. The formula witnessing failure of SPP for FO(n) in \cite{AvBNSub95} uses one $n$-ary relation symbol and it is equivalent to a universal FO(n+1)-formula.  Our result in Section 1  shows  that one binary relation symbol suffices for $n=3$.   Problems 4.3 and 4.2 in \cite{AvBNSub95} (restated in \cite[(2),(3) on p.701]{AvBNIGPL95} and \cite[(2) on p.241]{AvBNJPhL98}) ask whether just one binary relation suffices for all $n\ge 3$ and whether there is a bound on how many extra variables are needed for a universal equivalent of the witnessing formula (which always exists by the {\L}os-Tarski theorem for the language of FO as a whole). 
Proposition 23 in \cite{RW95} states that for all $n$, there is a FO(3) formula using only one binary relation symbol that is preserved by submodels but is not equivalent to any universal  FO(3+n)-formula, even on the finite models. This strong theorem  answers Problems 4.2,  4.3 in \cite{AvBNSub95}. Note: The examples used in \cite{RW95} supersede the  example presented in Section 1 of this note,  but they are built on a different idea, and we felt that the original example was still worth stating.

These are all results in the universe of all models. A classical result says that first-order logic over only finite models lacks the SPP, see \cite{Tait}. 

\vspace{1ex}

The results reviewed so far are for first-order logic with equality. What happens without equality in the language? At the end of the previous section, it is shown that SPP does not hold for equality-free FO(3). The same argument applied to the witness formulas in \cite[Proposition 23]{RW95} shows that equality-free FO(n) with $n\ge 3$ does not have SPP. The witness formula in \cite{GR99}, however, seems to use equality in an essential way. This leads to the following question which seems of some interest. Does equality-free two-variable first-order logic have the submodel preservation property? 
	
\vspace{1ex} 

Next, we consider some other types of fragments of first-order logic. Both the basic modal logic as well as the related first-order logic with generalized modal-style semantics have the SPP (see \cite[Theorem 2.10]{AvBNIGPL95}, \cite{AvBNJPhL98} for modal logic and \cite[p.712]{AvBNIGPL95}  for generalized semantics). 

We turn to extended modal languages. Fragments 1 and 2 in \cite{AvBNIGPL95} have the SPP, where Fragment 2 is the Guarded Fragment as now usually understood. The latter result is in \cite[Theorem 4.4.1]{AvBNJPhL98} while SPP also holds for the $n$-variable Guarded Fragment,  \cite[Theorem 4.4.2]{AvBNJPhL98}. Fragment 3 in \cite{AvBNIGPL95} is close to what is often called the Bounded Fragment or the basic hybrid logic with binder in the literature on hybrid logics, for which SPP seems still open, though a positive answer seems plausible given the techniques in \cite{BtC}, \cite{AtC07}.

\vspace{1ex}

A preservation property quite similar to, and often treated together with, the {\L}os-Tarski preservation property is the homomorphism preservation property (HPP). It says that a formula is preserved under homomorphisms if and only if it is equivalent to a positive existential formula  built from atomic ones by the use of conjunction, disjunction and existential quantifiers. Rosen \cite[Proposition 23, p.41]{R95} proves that HPP   holds for FO(2), and it was a long-standing open problem whether it holds for FO(n) for $n\ge 3$ (see, e.g., \cite[p.325]{GR99}). Bova and Chen \cite[Corollary 24]{BC19} settle this problem: unlike the {\L}os-Tarski preservation property, HPP does hold for FO(n) for all $n\ge 2$.

\vspace{1ex}

Finally, there is a literature on generalized versions of preservation theorems. For instance,  \cite{RW95} considers versions where the required syntactic equivalents may be in a specified extension of the base language. \cite{BB99} considers changes in the formulation of classical  system properties like Craig interpolation which are equivalent to the usual versions for first-order logic, but which generalize to infinitary logics where standard metatheorems fail. It is also worth noting that suitably generalized interpolation theorems, such as the powerful version in \cite{Otto2000}, imply Los-Tarski-style preservation theorems. 

\vspace{4ex}

\noindent {\bf Acknowledgment}\, We thank Balder ten Cate and Scott Weinstein for their helpful references and comments.

\bigskip\bigskip\bigskip\bigskip\bigskip\bigskip\bigskip

\noindent 
Andr\'eka, H.\ and N\'emeti, I.\\
Alfr\'ed R\'enyi Institute of Mathematics\\
Budapest, Re\'altanoda st.\ 13-15, H-1053 Hungary\\
andreka.hajnal@renyi.hu, nemeti.istvan@renyi.hu
\bigskip

\noindent
van Benthem, J.\\
Institute for Logic, Language and Computation, 
 University of Amsterdam\\
 P.O. Box 94242, 1090 GE Amsterdam, The Netherlands\\
 j.vanbenthem@uva.nl\\
Department of Philosophy, Stanford University, Stanford, CA 94305, USA\\
Department of Philosophy, Tsinghua University, Beijing, China

\end{document}